\RequirePackage{fix-cm}
\documentclass[smallextended]{svjour3}       
\smartqed  
\usepackage{graphicx}
\usepackage{graphics}
\usepackage{epsfig}
\usepackage{amsmath}
\usepackage{amssymb}

\newtheorem{prop}{Proposition}
\newcommand{\finish}{\hfill$\Box$\vspace{0.2cm}}
\newcommand{\rg}{\rightarrow}
\newcommand{\prf}{\noindent{\bf Proof:\ }}

\newcommand{\E}{{\rm I \!E}}
\newcommand{\p}{{\rm I \!P}}
\begin{document}

\title{Distribution of Maximum Loss for Fractional Brownian Motion
}

\titlerunning{Maximum Loss for fBm}        

\author{Mine Caglar         \and
        Ceren Vardar 
}


\institute{M. Caglar \at
              Department of Mathematics\\
              Koc University\\ Sariyer, Istanbul
              Tel.: +90-212-3381315\\
              Fax: +90-212-3381559\\
              \email{mcaglar@ku.edu.tr}           
           \and
           Ceren Vardar \at
              Department of Mathematics\\
              TOBB Economy and Technology University\\ Sogutozu,
              Ankara
              Tel.: +90-312-2924146\\
              Fax.: +90-312-2924324\\
              \email{cvardar@etu.edu.tr}
}

\date{Received: date / Accepted: date}

\maketitle

\begin{abstract}
In finance, the price of a volatile asset can be modeled using
fractional Brownian motion (fBm) with Hurst parameter
$H>\frac{1}{2}.$ The Black-Scholes model for the values of returns
of an asset using fBm is given as,
\[ Y_t=Y_0 \exp{((r+\mu)t+\sigma B_t^H)}\, ,~~~~t\geq 0\]
 where $Y_0$ is the initial value, $r$ is constant
interest rate, $\mu$ is constant drift and $\sigma$ is constant
diffusion coefficient of fBm, which is denoted by $(B_t^H)$ where $t
\geq 0.$ Black-Scholes model can be constructed with some Markov
processes such as Brownian motion. The advantage of modeling with
fBm to Markov proccesses is its capability of exposing the
dependence between returns. The real life data for a volatile asset
display long-range dependence property. For this reason, using fBm
is a more realistic model compared to Markov processes. Investors
would be interested in any kind of information on the risk in order
to manage it or hedge it. The maximum possible loss is one way to
measure highest possible risk. Therefore, it is an important
variable for investors. In our study, we give some theoretical
bounds on the distribution of maximum possible loss of fBm. We
provide both asymptotical and strong estimates for the tail
probability of maximum loss of standard fBm and fBm with drift and
diffusion coefficients. In the investment point of view, these
results explain, how large values of possible loss behave and its
bounds.
\end{abstract}

\section{Introduction}
\label{intro}

In finance, the price of one share of the risky asset, is modeled
using fractional Brownian motion (fBm) with Hurst parameter $H \in
[1/2,1)$ in order to display long-range dependence. Hence, our
results are on fractional Brownian motion with Hurst parameter $H
\in [1/2,1).$

 There are several stochastic
integral representations that have been developed for  fBm. For
example,
\begin{eqnarray}
B^{H}_{t}&=&\frac{1}{\Gamma(H+1/2)}\int_{\mathbb{R}}^{}((t-s)^{H-1/2}_{+}-(-s)^{H-1/2}_{+})dW_{s}\nonumber\\
&=&\frac{1}{\Gamma(H+1/2)}(\int_{-\infty}^{0}((t-s)^{H-1/2}-(-s)^{H-1/2})dW_{s}\nonumber\\
&+&\int_{0}^{t}(t-s)^{H-1/2}dW_{s})
\end{eqnarray}
is a fBm with Hurst parameter $H\in (0,1),$where $W$ is a
Wiener process.

Let $H$ be a constant in the interval $(0,1)$. A (standard) fBm $\{B_t^H: t \geq 0\}$ with Hurst parameter $H$ is a continuous and centered Gaussian
process with covariance function
\begin{equation}E[B_t^H B_s^H]=\frac{1}{2}(t^{2H}+s^{2H}-|t-s|^{2H})\end{equation}
and  $B_0^H=0$. It follows that $B^H$ has stationary increments, that is $B^H_{t+s}-B^H_s$ has the same law as $B_t^H,$ for $s,t\geq 0$.

  For $H=1/2$, the process
$(B_t^H)_{t\geq 0}$ corresponds to a standard Brownian motion, in
which the increments are independent. By  definition, the covariance
between the increments $B^H(t+h)-B^H(t)$  and $B^H(s+h)-B^H(t)$ with
$s+h\leq t$ and $t-s=nh$ is
\[\rho_H(n)=\frac{1}{2}h^{2H}[(n+1)^{2H}+(n-1)^{2H}-2n^{2H}].\]
We observe that two increments of the form $B^H(t+h)-B^H(t)$ and
$B^H(t+2h)-B^H(t+h)$ are positively correlated for $H>1/2$, and they
are negatively correlated for $H<1/2$.

Since the covariance function of fBm is homogeneous of order
$2H,$\cite{Biagini} fBm possesses the self-similarity property, that
is, for any constant $c>0,$
\[(B_{ct}^H)_{t\geq0}\stackrel{law}{=}(c^HB_t^H)_{t\geq0}.\]

 The aim of the present paper is to find the distribution of
maximum possible loss of fBm, which is one way of measuring
investors risk. In order to obtain results related to this
distribution, we are also interested in the supremum, the infemum, and
the range variables of fBm.

Let $(\Omega ,\Im ,P)$ be a probability space and let $B^H$ be  fractional Brownian motion with Hurst parameter $H \in
[1/2, 1).$ We introduce the following notation.

\begin{itemize}
\item  Let $I_t^H\mathop {: = }\limits^{} \mathop {\inf }\limits_{0 \le
v \le t} B_v^{}$ denote the \emph{infemum} of fractional Brownian motion up to
time $t.$

\item Let $S_t^H\mathop {: = }\limits^{} \mathop {\sup } \limits_{0 \le
v \le t} B_v^{}$ denote the \emph{supremum} of fractional Brownian motion
up to time $t.$

\item Let  $R_t^H = S_t^H - I_t^H$, called the \emph{range} of fractional Brownian motion up to time $t$.

\item The \emph{maximum loss} of fractional Brownian motion before time $t$ is defined as
    \[M_t^{H, - }: = \mathop {\sup }\limits_{0 \le u \le v \le
t} (B_u^{} - B_v^{}) = \mathop {\sup }\limits_{0 \le v \le t}
(\mathop {\sup }\limits_{0 \le u \le v} ({B_u} - {B_v}))\]
\end{itemize}

 More generally, the above definitions can be repeated for a non-standard fBm. Let $\mu$ be drift parameter taking real values other than $0$ and similarly $\sigma >0$ be a real valued diffusion
coefficient other than $1$ for fBm with drift defined as ${Y_t} := \mu t + \sigma {B_t}$.

\section{Bounds on the distribution of maximum loss}

In this section we introduce some new bounds on the expected value
of maximum loss of fBm and on the distribution of maximum loss of
fBm. These results are already useful for investors as they are.
Also, in later sections they will be useful for finding the
asymptotic distribution of maximum loss.

\begin{theorem}{\label{thm:1}} For fBm up to time $a$ with Hurst parameter
$H>\frac{1}{2},$ and for $y>0$, we have
 $$\frac{{\sqrt 2 {a^H}}}{{2\sqrt \pi
}} \le E(M_a^{H, - }) \le \frac{{2\sqrt 2 {a^H}}}{{\sqrt \pi  }}
$$ and $$P(M_a^{H, - }>y)<P(R_a^H \ge y) \le \frac{{2\sqrt 2
{a^H}}}{{y\sqrt \pi  }}$$
\end{theorem}
\prf $E(S_a^H) \le \frac{{\sqrt 2 }}{{\sqrt \pi  }} * {a^H}$ is
known \cite{Ceren} and by the Markov's inequality an upper bound for
the distribution of the supremum is found, that is for $x>0,$
$P(S_a^H>x)\leq\frac{\sqrt{2}a^H}{x\sqrt{\pi}}.$ Then, by the
symmetry property of centered Gaussian processes one can show that
$E(I_a^H) \ge - \frac{{\sqrt 2 }}{{\sqrt \pi  }} * {a^H}$

Combining the results given above we find an upper bound for the
expected value of range, $R_t^H,$ that is $E(R_a^H) \le
\frac{{2\sqrt 2 }}{{\sqrt \pi  }} * {a^H}.$ Furthermore, by Markov's
inequality we see that for $y>0,$ $P(R_a^H \ge y) \le \frac{{2\sqrt
2 {a^H}}}{{y\sqrt \pi}}.$

Clearly one can see that, $$I_t^H\mathop {: = }\limits^{} - \mathop
{\inf }\limits_{0 \le v \le t} X_v^{} \le \mathop {\sup }\limits_{u
\le v \le t} \mathop {\sup }\limits_{0 \le u \le v} (X_u^H - X_v^H)
= M_t^{H, - }\le R_t^H$$ holds for all $t.$

Hence \begin{equation}{\label{eq:1}} - E(I_a^H) \le E(M_a^{H, - })
\le E(R_a^H) \le \frac{{2\sqrt 2 {a^H}}}{{\sqrt \pi }}\end{equation}
is obtained, and by Markov's inequality we get
$$P(M_a^{H, - }>y)<P(R_a^H \ge y) \le \frac{{2\sqrt 2
{a^H}}}{{y\sqrt \pi  }}$$


We have also noticed that,
   $\frac{{\sqrt 2 {a^H}}}{{2\sqrt \pi  }} \le E(S_a^H) \le \frac{{\sqrt 2 {a^H}}}{{\sqrt \pi
   }}$ \cite{Norros}. These bounds are obtained using $H = 1$ and $H=1/2$ in Sudakov-Fernique
   inequality \cite[Theorem II.2.9]{Adler}.

    Furthermore, because $E(S_a^H)$ equals $- E(I_a^H)$ using Equation (\ref{eq:1})
    we obtain,
$$\frac{{\sqrt 2 {a^H}}}{{2\sqrt \pi  }} \le E(M_a^{H, - }) \le
\frac{{2\sqrt 2 {a^H}}}{{\sqrt \pi  }}.$$ \finish

\section{Asymptotical distribution of maximum loss of fractional Brownian motion}

The theorem in this section gives the asymptotic result for the probability
distribution of maximum loss of fBm. We first start with scaling property of loss process
which we denote as ${X_v}: = \mathop {\sup }\limits_{0 \le u \le v}
({B_u} - {B_v})$, $v>0.$ Our approach is similar to the large
deviations technique used for queeing systems modeled by fractional
Brownian motion, \cite{Duffield} and \cite{Norros}.

We define the loss process $X$ by $X_v=\mathop {\sup }\limits_{0 \le
u \le v} ({B_u} - {B_v})$.

\begin{prop}{\label{prop:1}}  The loss process $X$ is self-similar and each $X_v$ has the same
distribution as the supremum $S_v$ for every $v\geq 0$.
 \end{prop}
\prf The self-similarity of fBm corresponds to $\{ {B_{au}}:u \ge
0\} \mathop  = \limits^d \{ {a^H}{B_u}:u \ge 0\} $ for every $a>0$.
It follows that
\[\{ {B_{au}} - {B_{av}}:0 \le u \le v,v \ge 0\} \stackrel{d}{=}
\{ {a^H}({B_u} - {B_v}):0 \le u \le v,v \ge 0\} \; .
\]
Therefore, we get
\[\{ {X_{av}}:v \ge 0\}  \stackrel{d}{=}  \{ {a^H}{X_v}:v \ge 0\} \]
by the definition $X_v=\mathop {\sup }\limits_{0 \le u \le v} ({B_u}
- {B_v})$.

 On the other hand, since fractional Brownian motion has stationary increments, the
collections
 $\{ {B_u} - {B_v}:0 \le u \le v\} $ and  $\{ -{B_{v - u}}:0 \le
u \le v\} $ have the same probability law for fixed $v$. Both  are 0
mean Gaussian processes with covariance function
\[r(u,u') = 1/2\,[\;|v - u{|^{2H}} + |v - u'{|^{2H}} - |u - u'{|^{2H}}] \: .\]
Since the supremum of the two collections will also have the same
distribution and $\{  - {B_{v - u}}:0 \le u \le v\}  \mathop  =
\limits^d \{ {B_u}:0 \le u \le v\} $,  we get ${X_v}\mathop  =
\limits^d S_v^H$.
 \finish

Now, let us use $\Phi$ to denote cumulative distribution function of
Normal distribution and $\bar \Phi $ for its complement distribution

\begin{prop}{\label{prop:2}} For all $x \in {R_ + }$, we have
$$P(M_t^{H, - } > x) \ge \bar \Phi (x/{t^H}).$$
\end{prop}
 \prf As a result of Proposition \ref{prop:1}, for $0 \le v \le t,$
 we see that,
\begin{equation}P({X_v} > x) = P(\mathop {\sup }\limits_{0 \le u \le v} {B_u} >
x) \ge \mathop {\sup }\limits_{0 \le u \le v} P({B_u} > x) = \mathop
{\sup }\limits_{0 \le u \le v} \bar \Phi (x/{u^H}) = \bar \Phi
(x/{v^H}).\end{equation} Similarly, we observe
\begin{equation}P(M_t^{H, - }
> x) = P(\mathop {\sup }\limits_{0 \le v \le t} {X_v} > x) \ge
\mathop {\sup }\limits_{0 \le v \le t} P({X_v} > x) = \mathop {\sup
}\limits_{0 \le v \le t} \bar \Phi (x/{v^H}) = \bar \Phi
(x/{t^H})\end{equation}

\begin{theorem}{\label{thm:2}}For the maximum loss $M_t^{H, - }$ of fBm  and  $x>0,$ we have
$$\mathop {\lim }\limits_{x \to \infty
} \frac{1}{{{x^2}}}\,\,\log P(M_t^{H, - } > x) =  -
\frac{1}{{2{t^{2H}}}}$$
\end{theorem}

\prf We start with finding a lower bound for limit infemum of the
logarithm of the distribution of maximum loss. Combining Proposition \ref{prop:2}
with $$\mathop {\lim }\limits_{x \to \infty } (1/{x^2})\log \bar
\Phi (x/{t^H}) = - {(2{t^{2H}})^{ - 1}}$$(Adler, 1990, sf.42), we
obtain, \begin{equation}\mathop {\lim \inf }\limits_{x \to \infty }
\frac{1}{{{x^2}}}\log P(M_t^{H, - } > x) \ge \mathop {\lim
}\limits_{x \to \infty } \frac{1}{{{x^2}}}\log \bar \Phi (x/{t^H}) =
-\frac{1}{{2{t^{2H}}}}\end{equation}
We continue with finding an upper bound for the limit supremum.
Note that maximum loss $M_t^{H, - }$ has two dimensional index set
$T:=\{(u,v):0 \le u \le v \le t\},$ and it is a centered Gaussian
process. Based on the fact that $T$ is a separable metric space we use
Borel's inequality given in \cite[Theorem II.2.1]{Adler} for finding
the supremum of this process. By the continuity property of the
paths of fractional Brownian motion,
 ${B_u} - {B_v},$ is bounded on $T.$ And by Borel's inequality or directly from Theorem \ref{thm:1} we see that
$\eta : = E\,[\mathop {\sup }\limits_{0\le u \le v \le t} ({B_u} -
{B_v})]$ is finite.

As a result of Borel's inequality $$P(M_t^{H, - } > x)  = P(\mathop
{\sup }\limits_{0 \le u \le v \le t} ({B_u} - {B_v}) > x) \le 2{e^{
- \frac{1}{2}{{(x - \eta )}^2}/{t^{2H}}}}, \qquad x > \eta
$$ can be written. Here, $\mathop {\sup }\limits_{0 \le u \le v \le t}
{\rm{E}}{({B_u} - {B_v})^2} = {t^{2H}}$ is used.

This shows, \begin{equation}\mathop {\lim \sup }\limits_{x \to
\infty } \frac{1}{{{x^2}}} \log P(M_t^{H, - } > x) \le \mathop {\lim
}\limits_{x \to \infty } \left[ {\frac{{\log 2}}{{{x^2}}} -
\frac{1}{{{x^2}}} \frac{{{{(x - \eta )}^2}}}{{2{t^{2H}}}}} \right] =
- \frac{1}{{2{t^{2H}}}}\end{equation} which completes the proof.
\finish
\\
\\
In \cite{Marcus} it was shown that the result given in Theorem
\ref{thm:2} for centered, Gaussian sequence of random variables
which was mentioned in \cite[pg.43]{Adler}. However, in order for us
to use this result, it must be shown for general metric space $T,$
similar to the proof of Borel's inequality. In stead of showing this
generalization, we preferred giving the proof for $M_t^{H, - }.$ As
a conclusion, we see that the increase in the logarithm of
distribution of $M_t^{H, - }$ behaves same as the asymptotic
increase of distribution of ${B_t}.$ Here, the important property of
${B_t}$ is, it has the same distribution as the random variable with
the highest variation among the random variables $\{ {B_u} - {B_v}:0
\le u \le v \le t\}$ which is ${B_0} - {B_t} =  - {B_t}.$

\section{Asymptotical Estimate for the  Tail Probability}

In the next theorem, we generalize the previous result for fractional Brownian motion $Y$ with drift and diffusion coefficient. Recall
 ${Y_t} = \mu t + \sigma {B_t}$ for $\mu
\in R,\;\sigma  > 0.$

\begin{theorem} The maximum loss process defined by  $M_t^{H, - }: = \mathop
{\sup }\limits_{0 \le u \le v \le t} (Y_u^{} - Y_v^{})$ satisfies
\[\mathop {\lim }\limits_{x \to \infty } \frac{1}{{{x^2}}}\,\,\log P(M_t^{H, - } > x) =  - \frac{1}{{2\sigma^2{t^{2H}}}}\]
where $Y$ is a fractional Brownian motion with drift and diffusion
coefficients.
\end{theorem}

\prf For all $x \in {R_ + }$, the trivial lower bound is given by
\[P(M_t^{H, - } > x) = P(\mathop {\sup }\limits_{0 \le u \le v \le t} (Y_u^{} - Y_v^{}) > x) \ge \mathop {\sup }\limits_{0 \le u \le v \le t} P(Y_u^{} - Y_v^{} > x) \; .\]
Since $Y$ has stationary increments, we have
\[P(Y_u^{} - Y_v^{} > x) = P( - Y_{v - u}^{} > x) = \bar \Phi ((x + \mu (v - u))/(\sigma {(v - u)^H}))\; .
\]
Therefore, the following hold
\[\mathop {\sup }\limits_{0 \le u \le v \le t} P(Y_u^{} - Y_v^{} > x)
= \mathop {\sup }\limits_{0 \le u \le v \le t} \bar \Phi \left(
{\frac{{x + \mu (v - u)}}{{\sigma {{(v - u)}^H}}}} \right) = \mathop
{\sup }\limits_{0 \le v \le t} \bar \Phi \left(\frac{x + \mu
v}{\sigma {v^H}}\right)
\]
 Since $\Phi$ is a decreasing function, we find $v \in [0,t]$ that
 minimizes $f(v) = (x + \mu v)/(\sigma {v^H})$. The critical
 value of $f$ is $v* = \frac{{xH}}{{\mu (1 - H)}} $. We check if $v* \in
[0,t]$ or not as below
\begin{itemize}
\item If $\mu  > 0$ and $x < t\mu (1 - H)/H$, then the minimum value
of $f$ is obtained at $v* \in [0,t]$.
\item If $\mu  > 0$ and $x > t\mu (1 - H)/H$, then the minumum of $f$
is obtained at $t$.
\item If $\mu < 0$, then the minimum of  $f$ occurs at $t$.
\end{itemize}

As a result, we have $P(M_t^{H, - } > x) \ge \bar \Phi ((x + \mu
{\kern 1pt} t)/(\sigma {\kern 1pt} {t^H}))$ for large $x > 0$. We
get
\[\mathop {\lim \inf }\limits_{x \to \infty }  \frac{1}{{{x^2}}}\log P(M_t^{H, - } > x) \ge \mathop {\lim }\limits_{x \to \infty } \frac{1}{{{x^2}}}\log \bar \Phi (\frac{{x + \mu {\kern 1pt} t}}{{\sigma {\kern 1pt} {t^H}}}) =  - \frac{1}{{2{\sigma
^2}{t^{2H}}}} \; .
\]

On the other hand, we use Borel inequality which applies to a mean
zero Gaussian process in order to find the limit supremum
\cite[Theorem II.2.1]{Adler}. By definition, $Y_u^{} - Y_v^{} =
\sigma ({B_u} - {B_v}) + \mu (u - v)$ and
\begin{eqnarray*}
\mathop {\sup }\limits_{0 \le u \le v \le t} (Y_u^{} - Y_v^{}) > x &
\Leftrightarrow &  \exists (u,v) \in t\;, \sigma ({B_u} - {B_v}) +
\mu (u - v) > x \\
& \Leftrightarrow & \exists (u,v) \in t\;,
 \frac{{\sigma ({B_u} - {B_v})}}{{x - \mu (u -
v)}} > 1 \\ & \Leftrightarrow  & \mathop {\sup }\limits_{0 \le u \le
v \le t} \frac{{\sigma ({B_u} - {B_v})}}{{x - \mu (u - v)}} > 1
\end{eqnarray*}
where we assume $x > \mu t$. Making a change of variable $(u,v) \to
(u/x,v/x)$ and using the self-similarity of standard fractional
Brownian motion, we get
\begin{eqnarray*}
P(\mathop {\sup }\limits_{0 \le u \le v \le t} (Y_u^{} - Y_v^{}) >
x) & =& P\left( {\mathop {\sup }\limits_{0 \le u \le v \le t}
\frac{{\sigma ({B_u} - {B_v})}}{{x - \mu (u - v)}} > 1} \right) \\
& =& P\left( {\mathop {\sup }\limits_{0 \le u \le v \le \frac{t}{x}}
\,\,\frac{{\sigma ({B_{xu}} - {B_{xv}})}}{{x - \mu {\kern 1pt} x(u -
v)}} > 1} \right) \\
&= & P\left( {\mathop {\sup }\limits_{0 \le u \le v \le \frac{t}{x}}
\,\,\frac{{\sigma \,{x^H}({B_u} - {B_v})}}{{1 - \mu {\kern 1pt} (u -
v)}} > x} \right) \; .
\end{eqnarray*}
Since $G_{u,v}^{}: = \sigma \,({B_u} - {B_v})/(1 - \mu {\kern 1pt}
(u - v))$ is a zero-mean continuous Gaussian process,  Borel
inequality implies that for $x > {\eta _1}$
\begin{equation} P(M_t^{H, - } > x)  \label{ustsinir}
= P(\mathop {\sup }\limits_{0 \le u \le v \le t} (Y_u^{} - Y_v^{}) >
x) \le 2\exp \left[ { - \frac{1}{2}\frac{{{{\left( {{x^{1 - H}} -
{\eta _x}} \right)}^2}}}{\gamma }} \right]
\end{equation}
 where ${\eta _x} = E\,[\,\mathop {\sup
 }\limits_{0 \le u \le v \le t/x} G_{u,v}^{}]$, $\gamma  = \mathop {\sup }\limits_{0 \le u \le v \le t/x} E\,G_{u,v}^2$, and it is observed
 that ${\eta _x} \le {\eta _1}$ holds for large $x$, in particular $x
> t$. We have
\[\mathop {\sup }\limits_{0 \le u \le v \le t/x} G_{u,v}^{}\mathop  = \limits^d \mathop {\sup }\limits_{0 \le u \le v \le t/x} \frac{{  \sigma \,{B_{v - u}}}}{{1 + \mu {\kern 1pt} (v - u)}} = \mathop {\sup }\limits_{0 \le v \le t/x} \frac{{  \sigma \,{B_v}}}{{1 + \mu {\kern 1pt} v}}\]
and therefore
\begin{equation} \label{2}
0\le {\eta _x} \le \mathop {\sup }\limits_{0 \le v \le t/x} \frac{{  \sigma \,}}{{1 + \mu {\kern 1pt} v}}\,\,E(\mathop {\sup }\limits_{0 \le v \le t/x} {B_v})\; .
\end{equation}
Note that $ \sigma /(1 + \mu {\kern 1pt} v)$  is bounded from above
by $ \sigma$ for large $x$.  Since
\[\frac{{\sqrt 2 {{(t/x)}^H}}}{{2\sqrt \pi  }} \le E(\mathop {\sup }
\limits_{0 \le v \le t/x} {B_v}) \le \frac{{\sqrt 2
{{(t/x)}^H}}}{{\sqrt \pi  }}, \]
we take ${\eta _x} \propto {(t/x)^H}$ below in view of (\ref{2}). On the other hand,
 when $x$ is large enough, that is, when $t/x < H/(1 - \mu H)$,
\[\gamma  = \mathop {\sup }\limits_{0 \le u \le v \le t/x} E\,G_{u,v}^2 =
\mathop {\sup }\limits_{0 \le v \le t/x} \frac{{{\sigma ^2}
{v^{2H}}}}{{{{(1 + \mu v)}^2}}} = \frac{{{\sigma ^2}
{{(t/x)}^{2H}}}}{{{{(1 + \mu \,t/x)}^2}}}\] is found. As a result of
the bound (\ref{ustsinir}), we get
\[\mathop {\lim \sup }\limits_{x \to \infty } \frac{1}{{{x^2}}}\log P(M_t^{H, - } > x) \le \mathop {\lim }\limits_{x \to \infty } \left[ { - \frac{1}{{{x^2}}}\frac{{{{({x^{1 - H}} - {{(t/x)}^H})}^2}}}{{2{\sigma ^2}{{(t/x)}^{2H}}}}} \right] =  - \frac{1}{{2{\sigma ^2}{t^{2H}}}}\]
in view of $H \in (0,1)$. \finish

\section{Stronger Form of the Asymptotic Distribution}

In this section, we directly find the asymptotical form of the tail
probability using a characterization of \cite{Talagrand} for
Gaussian processes.

\begin{theorem} Asymptotically, we have
\[
\lim_{x\rg \infty} \frac{\p(M_t^{H-}>x)}{\Phi(x/t^{2H})}=1 \; .
\]
when $M_t$ is defined for standard fractional Brownian motion.
\end{theorem}
\prf
 The proof is based on two conditions given in \cite{Talagrand} that characterize the existence of the limiting distribution.
 We work with the mean zero Gaussian process
 $\{ {B_u} - {B_v}:0 \le u \le v \le
 t\}$. Let $\sigma^2_T=\sup_{(u,v)\in T}\E (B_{u}
 -B_{v})^2$, which yields $\sigma^2_T = t^{2H}$.
 The first condition is that there exists a unique $(u_0,v_0)$
 in $T$ such that $\E (B_{u_0} -B_{v_0})^2=\sigma_T^2$.
 This holds with $(u_0,v_0)=(0,t)$. For the second condition, the set
 \[
T_h=\{(u,v) \in T: \E(B_t(B_u -B_v)) \geq \sigma_T^2 -h^2 \}
 \]
is defined for $h>0$. In order to identify $T_h$, we find $ \E(B_t(B_u
-B_v))=\frac{1}{2} (v^{2H}-u^{2H}+(t-u)^{2H}-(t-v)^{2H})$.
Therefore, $(u,v) \in T_H$ satisfy
\begin{equation}
v^{2H}-u^{2H}+(t-u)^{2H}-(t-v)^{2H}\geq 2t^{2H}-2h^2\  \label{ineq}
\end{equation}
This implies that (\ref{ineq}) is satisfied by $(\bar{u},v) \in
T_H$, for  fixed $\bar{u}$. Now, for $\bar{u}> t/2$, we have
$(t-\bar{u})^{2H}-u^{2H}<0$, and $(t-v)^{2H}>0$ for all $v\in
[0,t]$. Then, from (\ref{ineq})
\[
v^{2H}\geq 2t^{2H}-2h^2 \geq t^{2H}-2h^2\; .
\]
On the other hand, for $\bar{u}\leq t/2$, we have $(t-\bar{u})^{2H}
\leq t^{2H}$. Therefore, we get
\[
v^{2H}+t^{2H} \geq v^{2H}-u^{2H}+(t-u)^{2H}-(t-v)^{2H}\geq
2t^{2H}-2h^2
\]
which again implies $v^{2H} \geq t^{2H}-2h^2$ for fixed $\bar{u}$
and $(\bar{u},v)\in T_h$. Since, $f(v)=v^{2H}$ is convex, we get
\[
v\geq t-Kh^2
\]
for some constant $K$ \cite[pg.309]{Talagrand}. Now, we consider
second condition in \cite{Talagrand} which requires
\begin{equation}  \label{ust}
\lim_{h\rg 0} h^{-1}\,  \E \sup_{({u},v)\in T_h} B_u -B_v+B_t =1
\end{equation}
In particular, we have
\[
\E \sup_{(\bar{u},v)\in T_h} B_{\bar{u}} -B_v+B_t \leq \E
\sup_{v\geq {\bar{u}}, \: t-v \leq Kh^2} B_{\bar{u}} -B_v+B_t
\]
Then, it follows that
\begin{equation} \label{alt}
\lim_{h\rg 0} h^{-1}\, \E \sup_{v\geq {\bar{u}}, \: t-v \leq Kh^2}
B_{\bar{u}} -B_v+B_t =  \lim_{h\rg 0} h^{-1}\, \E \sup_{ t-v \leq
Kh^2} B_{t-v}
\end{equation}
since fractional Brownian motion has stationary increments and $\E
B_{\bar{u}}=0$. On the right hand side of (\ref{alt}), the supremum
of fractional Brownian motion  $[0,Kh^2]$ is bounded by $\sqrt{2}K^H
h^{2H}/\sqrt{\pi}$ \cite{Ceren,Norros} and hence we get the limit in
(\ref{alt}) to be 0. Since $T$ is separable, a monotone convergence
argument extends the result for fixed $\bar{u}$ to all $T_h$ proving
(\ref{ust}) \cite[pg.47]{Adler}. \finish




\end{document}